\documentclass{article}
\usepackage[latin1]{inputenc}

\usepackage{abstract}
\usepackage{amsthm}
\newtheorem{theorem}{Theorem}
\newtheorem{proposition}{Proposition}

\usepackage{amssymb,amsmath}
\begin{document}
\title{Geometry of $CR$ submanifolds of maximal $CR$ dimension in complex space forms}
\author{Mirjana Milijevi\'{c}\\\\
Faculty of Architecture and Civil Engineering\\
University of Banja Luka\\
Stepe Stepanovi\'{c}a 77, 78000 Banja Luka\\
Bosnia and Herzegovina\\
E-mail: \texttt{$mirjana_{-}milijevic@yahoo.com$}}
\date{}
\maketitle
\begin{abstract}
On real hypersurfaces in complex space forms many results are proven. In this paper we generalize some results concerning extrinsic geometry of real hypersurfaces, to $CR$ submanifolds of maximal $CR$ dimension in complex space forms.
\end{abstract}
\section{Introduction}
Let $\overline{\textbf{M}}$ be an $(n + p)$-dimensional complex space form, i.e. a Kaehler manifold of constant holomorphic sectional curvature $4c$, endowed with metric $\overline{g}$. Let $\textbf{M}$ be an $n$-dimensional real submanifold of $\overline{\textbf{M}}$ and $J$ be the complex structure of $\overline{\textbf{M}}$. For a tangent space $T_{x}(\textbf{M})$ of $\textbf{M}$ at $x$,\:we put $H_{x}(\textbf{M})=JT_{x}(\textbf{M})\cap T_{x}(\textbf{M})$. Then,\:$H_{x}(\textbf{M})$ is the maximal complex subspace of $T_{x}(\textbf{M})$ and is called the holomorphic tangent space to $\textbf{M}$ at $x$. If the complex dimension $dim_{\textbf{C}}H_{x}(\textbf{M})$ is constant over $\textbf{M}$,\:$\textbf{M}$ is called a Cauchy-Riemann submanifold or briefly a $CR$ submanifold and the constant $dim_{\textbf{C}}H_{x}(\textbf{M})$ is called the $CR$ dimension of $\textbf{M}$. If,\:for any $x\in \textbf{M}$,\:$H_{x}(\textbf{M})$ satisfies $dim_{\textbf{C}}H_{x}(\textbf{M})=\frac{n-1}{2}$,\:then $\textbf{M}$ is called a $CR$ submanifold of maximal $CR$ dimension. It follows that there exists a unit vector field $\xi$ normal to $\textbf{M}$ such that $JT_{x}(\textbf{M})\subset T_{x}(\textbf{M})\oplus span\{\xi_{x}\}$,\:for any $x\in \textbf{M}$.

A real hypersurface is a typical example of a $CR$ submanifold of maximal $CR$ dimension. The study of real hypersurfaces in complex space forms is a classical topic in differential geometry and the generalization of some results which are valid for real hypersurfaces to $CR$ submanifolds of maximal $CR$ dimension may be expected.

For instance, nonexistence of real hypersurfaces with the parallel shape operator ([1], [2]) and real hypersurfaces with the second fundamental form satisfying $h(JX,Y)-Jh(X,Y)=0$ ([3]), in nonflat complex space forms, is proven.

In this paper we study the conditions that the shape operator of the distinguished vector field $\xi$ is parallel and that the second fundamental form satisfies $h(JX,Y)-Jh(X,Y)=0$, on $CR$ submanifolds of maximal $CR$ dimension in complex space forms.

The author wishes to express her gratitude to Professor Mirjana Djori\'{c} for her useful advice.

\section{$CR$ submanifolds of maximal $CR$ dimension of a complex space form}
Let $\overline{\textbf{M}}$ be an $(n + p)$-dimensional complex space form with Kaehler structure $(J,\overline{g})$ and of constant holomorphic sectional curvature $4c$. Let $\textbf{M}$ be an $n$-dimensional $CR$ submanifold of maximal $CR$ dimension in $\overline{\textbf{M}}$ and
$\iota :\textbf{M}\to {\overline{\textbf{M}}}$  immersion. Also,\:we denote by $\iota$ the differential of the immersion. The Riemannian metric $g$ of $\textbf{M}$ is induced from the Riemannian metric $\overline{g}$ of $\overline{\textbf{M}}$ in such a way that
$g(X,Y)=\overline{g}(\iota X,\iota Y)$,\:where $X,\:Y \in T(\textbf{M})$. We denote by $T(\textbf{M})$ and $T^{\bot}(\textbf{M})$ the tangent bundle and the normal bundle of $\textbf{M}$,\:respectively.

On $\overline{\textbf{M}}$ we have the following decomposition into tangential and normal components:
\begin{align}\label{eq:m1}
J\iota X=\iota FX+u(X)\xi,\:\: X\in T(\textbf{M}).
\end{align}
Here $F$ is a skew-symmetric endomorphism acting on $T(\textbf{M})$ and $u$ in one-form on $\textbf{M}$.

Since $T_{1}^{\bot}(\textbf{M})=\{\eta\in T^{\bot}(\textbf{M})|\overline{g}(\eta,\xi)=0\}$ is $J$-invariant,\:from now on we will denote the orthonormal basis of $T^{\bot}(\textbf{M})$ by $\xi,\xi_{1},\cdots,\xi_{q},\xi_{1^{*}},\cdots,\xi_{q^{*}}$,\:where $\xi_{a^{*}}=J\xi_{a}$ and $q=\frac{p-1}{2}$. Also,\:$J\xi$ is the vector field tangent to $\textbf{M}$ and we write

\begin{align}\label{eq:m2}
J\xi=-\iota U.
\end{align}
Furthermore,\:using (\ref{eq:m1}),\:(\ref{eq:m2}) and the Hermitian property of $J$ implies
\begin{align}\label{eq:80}
F^{2}X=-X+u(X)U,
\end{align}
\begin{align}\label{eq:m5}
FU=0,
\end{align}
\begin{align}\label{eq:m6}
g(X,U)=u(X).
\end{align}
Next,\:we denote by $\overline{\nabla}$ and $\nabla$ the Riemannian connection of $\overline{\textbf{M}}$ and $\textbf{M}$,\:respective-\-ly, and by $D$ the normal connection induced from $\overline{\nabla}$ in the normal bundle of $\textbf{M}$. They are related by the following Gauss equation
\begin{align}\label{eq:m7}
\overline{\nabla}_{\iota X}\iota Y=\iota \nabla_{X}Y+h(X,Y),
\end{align}
where $h$ denotes the second fundamental form,\:and by Weingarten equations
\begin{align}\label{eq:m8}
\overline{\nabla}_{\iota X}\xi&=-\iota AX+D_{X}\xi\\
\notag
&=-\iota AX+\sum_{a=1}^{q}\{s_{a}(X)\xi_{a}+s_{a^{*}}(X)\xi_{a^{*}}\},
\end{align}
\begin{align}\label{eq:m9}
\overline{\nabla}_{\iota X}\xi_{a}&=-\iota A_{a}X+D_{X}\xi_{a}=-\iota A_{a}X-s_{a}(X)\xi\\
\notag
&+\sum_{b=1}^{q}\{s_{ab}(X)\xi_{b}+s_{ab^{*}}(X)\xi_{b^{*}}\},
\end{align}
\begin{align}\label{eq:m10}
\overline{\nabla}_{\iota X}\xi_{a^{*}}&=-\iota A_{a^{*}}X+D_{X}\xi_{a^{*}}=-\iota A_{a^{*}}X-s_{a^{*}}(X)\xi\\
\notag
&+\sum_{b=1}^{q}\{s_{a^{*}b}(X)\xi_{b}+s_{a^{*}b^{*}}(X)\xi_{b^{*}}\},
\end{align}
where the $s$'s are the coefficients of the normal connection $D$ and $A$,\:$A_{a}$,\:$A_{a^{*}}$;\:$a=1,\cdots,q$,\:are the shape operators corresponding to the normals $\xi$,\:$\xi_{a}$,\:$\xi_{a^{*}}$,\:respecti-\-vely. They are related to the second fundamental form by
\begin{align}\label{eq:m11}
h(X,Y)&=g(AX,Y)\xi\\
\notag
&+\sum_{a=1}^{q}\{g(A_{a}X,Y)\xi_{a}+g(A_{a^{*}}X,Y)\xi_{a^{*}}\}.
\end{align}
Since the ambient manifold is a Kaehler manifold,\:using (\ref{eq:m1}),\:(\ref{eq:m2}),\:(\ref{eq:m9}) and (\ref{eq:m10}),\:it follows that
\begin{align}\label{eq:m12}
A_{a^{*}}X=FA_{a}X-s_{a}(X)U,
\end{align}
\begin{align}\label{eq:m13}
A_{a}X=-FA_{a^{*}}X+s_{a^{*}}(X)U,
\end{align}
\begin{align}\label{eq:m14}
s_{a^{*}}(X)=u(A_{a}X),
\end{align}
\begin{align}\label{eq:m15}
s_{a}(X)=-u(A_{a^{*}}X),
\end{align}
for all $X,\:Y$ tangent to $\textbf{M}$ and $a=1,\cdots,q$.\\
Moreover,\:since $F$ is skew-symmetric and $A_{a}$ and $A_{a^{*}}$;\:$a=1,\cdots,q$,\:are symmetric,\:(\ref{eq:m12}) and (\ref{eq:m13}) imply
\begin{align}\label{eq:m16}
g((A_{a}F+FA_{a})X,Y)=u(Y)s_{a}(X)-u(X)s_{a}(Y),
\end{align}
\begin{align}\label{eq:m17}
g((A_{a^{*}}F+FA_{a^{*}})X,Y)=u(Y)s_{a^{*}}(X)-u(X)s_{a^{*}}(Y),
\end{align}
for all $a=1,\cdots,q$.\\
Finally,\: the Codazzi equation for the distinguished vector field $\xi$ becomes
\begin{align}\label{eq:m21}
&(\nabla_{X}A)Y-(\nabla_{Y}A)X=c\{u(X)FY-u(Y)FX-2g(FX,Y)U\}\\
\notag
&+\sum_{a=1}^{q}\{s_{a}(X)A_{a}Y-s_{a}(Y)A_{a}X\}+\sum_{a=1}^{q}\{s_{a^{*}}(X)A_{a^{*}}Y-s_{a^{*}}(Y)A_{a^{*}}X\},
\end{align}
for all $X,\:Y$ tangent to $\textbf{M}$.

\section{Shape operator $A$ is parallel}
Here, we will give one well known result about hypersurfaces with the parallel shape operator.
\begin{theorem}\label{T1}[1], [2]
Let $\textbf{M}$ be an $n$-dimensional,\:where $n\geq 3$,\:hypersurface in a complex space form of constant holomorphic sectional curvature $4c \neq 0$. Then the shape operator $A$ of $\textbf{M}$ cannot be parallel.
\end{theorem}
We will study the same condition on $CR$ submanifolds of maximal $CR$ dimension in complex space forms. Therefore,\:we have the next two theorems.
\begin{theorem}\label{T2}
 Let $\textbf{M}$ be an $n$-dimensional $CR$ submanifold of maximal $CR$ dimension in an $(n+p)$-dimensional complex space form $(\overline{\textbf{M}},J,\overline{g})$, where $n\geq 3$ and the constant holomorphic sectional curvature of $\overline{\textbf{M}}$ equals $4c$.  Let the distinguished vector field $\xi$ be parallel with respect to the normal connection $D$ and $A$ be the shape operator of $\xi$. If $\nabla A=0$ on $\textbf{M}$,\:then $\overline{\textbf{M}}$ is an Euclidean space.
\end{theorem}
\flushleft \emph{Proof.}
Putting $Y=U$ in the Codazzi equation (\ref{eq:m21}),\:we get
\begin{align}
\notag
&(\nabla_{X}A)U-(\nabla_{U}A)X=-cFX +\sum_{a=1}^{q}\{s_{a}(X)A_{a}U-s_{a}(U)A_{a}X\}+\\
\notag
&\sum_{a=1}^{q}\{s_{a^{*}}(X)A_{a^{*}}U-s_{a^{*}}(U)A_{a^{*}}X\}.
\end{align}
From the assumption of the Theorem 2 and the last equation we get
\begin{align}\label{eq:m22}
&c FX=0.
\end{align}
From the equation (\ref{eq:m22}) we conclude that $c=0$.\:$\Box$
\begin{theorem}\label{T3}
Let $\textbf{M}$ be an $n$-dimensional $CR$ submanifold of maximal $CR$ dimension in an $(n+p)$-dimensional complex space form
$(\overline{\textbf{M}},J,\overline{g})$, where $n\geq 3$ and the constant holomorphic sectional curvature of $\overline{\textbf{M}}$ equals $4c$.
Let $p < n$ and $A$ be the shape operator of the distinguished vector field $\xi$. If $\nabla A=0$ on $\textbf{M}$, then $\overline{\textbf{M}}$ is an Euclidean space.
\end{theorem}
\flushleft \emph{Proof.}
 After putting $Y=U$ in (\ref{eq:m21}) and using the assumption of the Theorem 3, we get
 \begin{align}\label{eq:m81}
\notag
&-cFX+\sum_{a=1}^{q}\{s_{a}(X)A_{a}U-s_{a}(U)A_{a}X\}+\\
&\sum_{a=1}^{q}\{s_{a^{*}}(X)A_{a^{*}}U-s_{a^{*}}(U)A_{a^{*}}X\}=0.
\end{align}
 Multiplying the equation (\ref{eq:m81}) with an arbitrary vector field $Y\in T(\textbf{M})$,\:we get
\begin{align}\label{eq:m60}
&-cg(FX,Y)+\sum_{a=1}^{q}\{s_{a}(X)g(A_{a}U,Y)-s_{a}(U)g(A_{a}X,Y)\}+\\
\notag
&\sum_{a=1}^{q}\{s_{a^{*}}(X)g(A_{a^{*}}U,Y)-s_{a^{*}}(U)g(A_{a^{*}}X,Y)\}=0.
\end{align}
Interchanging $X$ and $Y$ in (\ref{eq:m60}) and subtracting (\ref{eq:m60}) and the resulting equation, we get
\begin{align}
\notag
&-2cg(FX,Y)+\sum_{a=1}^{q}\{s_{a}(X)g(A_{a}U,Y)+s_{a^{*}}(X)g(A_{a^{*}}U,Y)\}-\\
\notag
&\sum_{a=1}^{q}\{s_{a}(Y)g(A_{a}U,X)+s_{a^{*}}(Y)g(A_{a^{*}}U,X)\}=0.
\end{align}
Now,\:using (\ref{eq:m6}),\:(\ref{eq:m14}) and (\ref{eq:m15}),\:from the last equation it follows that
\begin{align}\label{eq:m61}
cFX=\sum_{a=1}^{q}\{s_{a}(X)A_{a}U+s_{a^{*}}(X)A_{a^{*}}U\}.
\end{align}
From (21) it follows that $FX$ is a linear combination of $A_{a}U$ and $A_{a^{*}}U$; $a=1,\cdots,q$.\\
Since every tangent vector $Y$ orthogonal to $U$ can be expressed as $Y=FX$, for $n-1>2q=p-1$, i.e. $n>p$, it follows that there exists a unit vector field $Y=FX$ which is orthogonal to $span\{A_{a}U,A_{a^{*}}U\}$; $a=1,\cdots,q$. For such $Y=FX$ it follows $s_{a}(Y)=0=s_{a^{*}}(Y)$; $a=1,\cdots,q$, using (13) and (14).\\
Consequently, using (21), we obtain
\begin{align}
\notag
cF^{2}X=\sum_{a=1}^{q}\{s_{a}(FX)A_{a}U+s_{a^{*}}(FX)A_{a^{*}}U\}.
\end{align}
Finally, using (3), we conclude $c=0$.\:$\Box$

\section{$CR$ submanifolds of maximal $CR$ dimension satisfying $h(JX,Y)=Jh(X,Y)$ }
On real hypersurfaces the next theorem is proven.
\begin{theorem}\label{T4}[3]
Let $\textbf{M}$ be an $n$-dimensional, $n\geq 3$, real hypersurface in a complex space form $(\overline{\textbf{M}},J,\overline{g})$. If the second fundamental form satisfies condition $h(JX,Y)=Jh(X,Y)$; $X, Y\in T(\textbf{M})$, then $\overline{\textbf{M}}$ is an Euclidean space.
 \end{theorem}
In the next theorem we will see if the same result is true on $CR$ submanifolds of maximal $CR$ dimension.
 \begin{theorem}\label{T5}
  Let $\textbf{M}$ be an $n$-dimensional, $n\geq 3$, $CR$ submanifold of maximal $CR$ dimension in an $(n+p)$-dimensional complex space form $(\overline{\textbf{M}},J,\overline{g})$ and the constant holomorphic sectional curvature of $\overline{\textbf{M}}$ equals $4c$. If the second fundamental form satisfies condition $h(JX,Y)=Jh(X,Y)$; $X, Y\in T(\textbf{M})$, then $\overline{\textbf{M}}$ is an Euclidean space.
  \end{theorem}
  \flushleft \emph{Proof.}
Using (\ref{eq:m11}),\:we have the next two equations
\begin{align}\label{eq:m23}
h(JX,Y)=g(AJX,Y)\xi + \sum_{a=1}^{q}\{g(A_{a}JX,Y)\xi_{a} + g(A_{a^{*}}JX,Y)\xi_{a^{*}}\}
\end{align}
and
\begin{align}\label{eq:m24}
Jh(X,Y)=-g(AX,Y)\iota U+\sum_{a=1}^{q}\{g(A_{a}X,Y)\xi_{a^{*}}-g(A_{a^{*}}X,Y)\xi_{a}\}.
\end{align}
From (\ref{eq:m23}),\:(\ref{eq:m24}) and the assumption of the Theorem 5,\:we have
\begin{align}\label{eq:m25}
&g(AX,Y)\iota U+g(JX,AY)\xi+\sum_{a=1}^{q}\{g(JX,A_{a}Y)+g(A_{a^{*}}X, Y)\}\xi_{a}+\\
\notag
&\sum_{a=1}^{q}\{g(JX,A_{a^{*}}Y)-g(A_{a}X,Y)\}\xi_{a^{*}}=0,
\end{align}
where we used the symmetry of the shape operators $A,\:A_{a},\:A_{a^{*}}$;\: $a=1,\cdots,q$.\\
From (\ref{eq:m1}) and (\ref{eq:m25}),\:we have
\begin{align}\label{eq:m26}
&g(AX,Y)\iota U+g(AFX,Y)\xi+\sum_{a=1}^{q}\{g(A_{a}FX,Y)+g(A_{a^{*}}X, Y)\}\xi_{a}+\\
\notag
&\sum_{a=1}^{q}\{g(A_{a^{*}}FX,Y)-g(A_{a}X,Y)\}\xi_{a^{*}}=0.
\end{align}
Because of the linear independence of the vectors
\begin{align}
\notag
\iota U,\:\xi,\:\xi_{a},\:\xi_{a^{*}};\:a=1,\cdots,q,
\end{align}
from (\ref{eq:m26}) we get the next equations
\begin{align}\label{eq:m27}
A=0,
\end{align}
\begin{align}\label{eq:m29}
A_{a}F=-A_{a^{*}};\:\: a=1,\cdots,q,
\end{align}
and
\begin{align}\label{eq:m30}
A_{a^{*}}F=A_{a};\:\: a=1,\cdots,q.
\end{align}
From the Codazzi equation (\ref{eq:m21}) and (\ref{eq:m27}),\:we conclude that
\begin{align}\label{eq:m53}
&0=c\{u(X)FY-u(Y)FX-2g(FX,Y)U\}\\
\notag
&+\sum_{a=1}^{q}\{s_{a}(X)A_{a}Y-s_{a}(Y)A_{a}X\}+\sum_{a=1}^{q}\{s_{a^{*}}(X)A_{a^{*}}Y-s_{a^{*}}(Y)A_{a^{*}}X\}.
\end{align}
Now,\:from the equations (\ref{eq:m12}) and (\ref{eq:m29}) we get
\begin{align}
\notag
(FA_{a}+A_{a}F)X=s_{a}(X)U;\:a=1,\cdots,q.
\end{align}
By scalar multiplication of the last equation with an arbitrary $Y\in T(\textbf{M})$, and using (\ref{eq:m16}), we get
\begin{align}\label{eq:m39}
s_{a}(Y)=0;\:a=1,\cdots,q.
\end{align}
From (\ref{eq:m13}) and (\ref{eq:m30}) we get
\begin{align}
\notag
(A_{a^{*}}F+FA_{a^{*}})X=s_{a^{*}}(X)U;\:a=1,\cdots,q.
\end{align}
By scalar multiplication of the last equation with an arbitrary $Y\in T(\textbf{M})$, and using (\ref{eq:m17}), we get
\begin{align}\label{eq:m40}
s_{a^{*}}(Y)=0;\:a=1,\cdots,q.
\end{align}
From (\ref{eq:m53}),\:(\ref{eq:m39}) and (\ref{eq:m40}),\:we get
\begin{align}
\notag
0=c\{g(X,U)FY-g(Y,U)FX-2g(FX,Y)U\}.
\end{align}
Multiplying the last equation with $U$,\:we get
\begin{align}\label{eq:m41}
0=-2cg(FX,Y).
\end{align}
From (\ref{eq:m41}) we conclude that $c=0$.\:$\Box$\\
Now, using (\ref{eq:m39}) and (\ref{eq:m40}),\:we get the next proposition.
\begin{proposition}
Let $\textbf{M}$ be an $n$-dimensional, $n\geq 3$, $CR$ submanifold of maximal $CR$ dimension in an $(n+p)$-dimensional complex space form $(\overline{\textbf{M}},J,\overline{g})$. If on $\textbf{M}$ the second fundamental form satisfies condition $h(JX,Y)=Jh(X,Y)$; $X, Y\in T(\textbf{M})$,\:then the distinguished vector field $\xi$ is parallel with respect to the normal connection $D$.
\end{proposition}

\end{document}